\documentclass[a4paper, 12pt]{amsart}
\usepackage{amsmath, amssymb, amscd, enumerate}

\newtheorem{lem}{Lemma}

\newtheorem{prop}{Proposition}
\newtheorem{conj}{Conjecture}

\newcommand{\pf}{{\noindent\bf Proof:\,}}
\makeatletter
\def\@eqnnum{(\thesection.\theequation)}
\makeatother

\begin{document}

\title{ On Nicolas criterion for the Riemann Hypothesis}

\author{YoungJu Choie}

\address{Department of Mathematics, Pohang Mathematics Institute,  POSTECH, Pohang, Korea}

\email{yjc@postech.ac.kr}

\author{Michel Planat }

\address{Institut FEMTO-ST, CNRS, 32 Avenue de
l'Observatoire, F-25044 Besan\c con, France }

\email{michel.planat@femto-st.fr }

\author{Patrick Sol\'{e} }
\address{Telecom ParisTech, 46 rue Barrault, 75634
Paris Cedex 13, France. }

\email{sole@enst.fr }

\subjclass[2000]{Primary 11F11,   the second author partially
supported by NRF 2009-0083-919 and NRF-2009-0094069}

\keywords{Nicolas inequality, Euler totient, Dedekind $\Psi$
function, Riemann Hypothesis, Primorial numbers}

\dedicatory{}

\begin{abstract}
Nicolas criterion for the Riemann Hypothesis is based on an
inequality that Euler totient function must satisfy at primorial
numbers. A natural approach to derive this inequality would be to
prove that a specific sequence related to that bound is strictly
decreasing. We show that, unfortunately, this latter fact would
contradict Cram\'er conjecture on gaps between consecutive primes.
An analogous situation holds when replacing Euler totient by
Dedekind $\Psi$ function.
\end{abstract}

\maketitle

\section{Introduction}
The Riemann Hypothesis (RH), which describes the non trivial
zeroes of Riemann $\zeta$ function has been qualified of Holy
Grail of Mathematics by several authors \cite{BCR,L}. There exist
many equivalent formulations in the literature \cite{C}. The one
of concern here is that of Nicolas \cite{N} that states that the
inequality

$$\frac{N_k}{\varphi(N_k)} >e^\gamma \log \log N_k, $$
where
\begin{itemize}
\item $\gamma \approx 0.577$ is the Euler Mascheroni constant,
\item $\varphi$ Euler totient function ,\item
$N_n=\prod_{k=1}^np_k$ the primorial of order $n,$
\end{itemize}
holds for all $k \ge 1$ if RH is true \cite[Th. 2 (a)]{N}.
Conversely, if RH is false, the inequality holds for infinitely
many $k,$ and is violated for infinitely many $k$ \cite[Th. 2
(b)]{N}. Thus, it is enough, to confirm RH, to prove this
inequality for $k$ large enough. In this note, we show that a
natural approach to this goal fails conditionally on a conjecture
arguably harder than RH, namely Cram\'er conjecture \cite{C}
$$p_{n+1}-p_n=O(\log^2 p_n).$$
Note that under RH, it can only be shown that \cite{C2}
$$p_{n+1}-p_n=O(\sqrt{ p_n}\log p_n).$$
See \cite{G} for a critical discussion of this conjecture. An
important ingredient of our proof is Littlewood oscillation
Theorem for Chebyshev $\theta$ function \cite[Th. 6.3]{I}. An
analogous situation holds when replacing Euler totient by Dedekind
$\Psi$ function, and replacing Nicolas criterion by \cite[Th.
2]{SP}.
\section{An intriguing sequence}
{\bf General conventions:}
\begin{enumerate}
\item We write $\log_2$ for  $\log \log,$ and $\log_3$ for $\log
\log_2$\item The formula $f=O(g)$ means that $\exists C>0,$ such
that $\vert f \vert \le C g.$  \item The formula $a_k \sim b_k$
means that $\forall \epsilon>0, \exists k_0,$
such that $b_k (1-\epsilon)\le a_k\le b_k (1+\epsilon),$ if $k> k_0.$\\
\end{enumerate}

We begin by an easy application of Mertens formula \cite[Th.
429]{HW}. For convenience define

$$ R(n)=\frac{n}{\varphi(n) \log_2 n }. $$
 Recall,
for future use, $\theta(x),$ Chebyshev's first summatory function:
$$\theta(x)= \sum_{p\le x }\log p.$$
 {\prop \label{mertens}  For $n$ going to $\infty$ we have $$\lim
R(N_n)=e^\gamma.$$}

\pf
 Put $x=p_n$ into Mertens
formula
$$\prod_{p\le x}(1-1/p)^{-1}\sim e^\gamma \log(x)$$ to obtain

$$R(N_n)\sim {e^\gamma} \log(p_n),$$
Now the Prime Number Theorem \cite[Th. 6, Th. 420]{HW} shows that
$x \sim \theta(x)$ for $x$ large. This shows that, taking $x=p_n$
we have
$$p_n\sim \theta(p_n)=\log(N_n).$$
The result follows. \qed

Define the sequence
$$u_n=R(N_n).$$

We have just shown that this sequence converges to $e^\gamma.$ But
Nicolas inequality is equivalent to saying that

$$u_n> e^\gamma.$$

So we observe

{\prop If $u_n$ is strictly decreasing for $n$ big enough then
Nicolas inequality is satisfied for $n$ big enough.}

\pf Assume $u_n >u_{n+1}$ for $n>n_0$ and that  Nicolas inequality
is violated for $N >n_0$ that is
$$u_n \le e^\gamma,$$

then for $n \ge N+1$ we have $u_{n+1} < u_n   \le e^\gamma.$ This
implies

$$\overline{\lim}\, u_n < e^\gamma,$$

contradicting Proposition \ref{mertens}.

\qed

We reduce the decreasing character of $u_n$ to a concrete
inequality between arithmetic functions.

{\prop The inequality $u_n>u_{n+1}$ is equivalent to
\begin{equation} \label{aritheq}
\log(1+\frac{\log p_{n+1} }{\theta(p_n)}
                                        )>\frac{\log
                                        \theta(p_{n+1})}{p_{n+1}}.
\end{equation}
}

\pf The inequality $u_n>u_{n+1}$ can be written as

$$\frac{N_{n}}{\varphi(N_{n})\log_2 N_n}  >\frac{N_{n+1}}{\varphi(N_{n+1})\log_2 N_{n+1}}.$$

 Note first that
$$ \frac{N_{n+1}}{\varphi(N_{n+1})} =\frac{1}{(1-1/p_{n+1})} \frac{N_{n}}{\varphi(N_{n})},$$

so that, after clearing denominators, $u_n>u_{n+1}$ is equivalent
to
$$ \log_2( N_{n+1})(1-1/p_{n+1})> \log_2 N_n, $$
or, distributing, to
$$ \log_2(N_{n+1} )- \log_2 N_n> \frac{\log_2 N_{n+1}}{p_{n+1}}. $$

 Now, to evaluate the LHS we write $N_{n+1}=N_n p_{n+1}$ so that
$$ \log_2(N_{n+1} )=\log_2( N_n p_{n+1})=\log ( \log N_n +\log p_{n+1})=\log_2 N_n+ \log(1+\frac{\log p_{n+1}}{\log N_n}).$$
to obtain
$$ \log(1+\frac{\log p_{n+1}}{\log N_n} )>\frac{\log_2 N_{n+1}}{p_{n+1}}.$$
 The result follows then upon letting $\log N_n=\theta(p_n).$
\qed

In fact, more could be true.

{\conj   Inequality (\ref{aritheq}) holds for all $n\ge 1.$ }

A heuristic motivation runs as follows $$\log(1+\frac{\log
p_{n+1}}{\theta(p_n)})\approx \frac{\log p_{n+1}}{\theta(p_n)}
\approx \frac{\log p_{n+1}}{p_n}.$$ Similarly

$$\frac{\log \theta(p_{n+1})}{p_{n+1}} \approx \frac{\log
p_{n+1}}{p_{n+1}}.$$ But, trivially

$$ \frac{\log p_{n+1}}{p_n}>\frac{\log
p_{n+1}}{p_{n+1}}.$$
 Numerical computations
confirm Conjecture 1 up to $n \le 10000.$ Unfortunately,
Proposition \ref{cram} provides a conditional disproof of this
conjecture.
\section{Background material}
We need an easy consequence of Littlewood oscillation theorem.

{\lem \label{Littlewood} There are infinitely many $n$ such that
 $$\theta(p_n)>k_n=p_n+C\sqrt{p_n} \log_3 p_n,$$ for some constant $C$
 independent of $n.$}

\pf By \cite[Th. 6.3]{I}, we know there are infinitely many values
of $x$ such that
$$\theta(x)>x+C\sqrt{x} \log_3 x.$$
Let $p_n$ be the largest prime $\le x.$  Thus
$$\theta(p_n)=\theta(x)>x+C\sqrt{x}>p_n+C\sqrt{p_n} \log_3 p_n.$$
\qed
\section{More on $u_n$}
Unfortunately, the sequence $u_n$ is not decreasing as the next
Proposition shows, conditionally on Cram\'er conjecture.

 {\prop \label{cram} The
inequality $u_n>u_{n+1}$ is violated for infinitely many $n$'s.}

 \pf By Lemma \ref{Littlewood} there are infinitely many $n$ such that $\theta(p_n)>k_n.$
For these $n$  the RHS of (\ref{aritheq}) is $> \frac{\log k_{n+1}}{p_{n+1}}>\frac{\log k_n}{p_{n+1}}.$\\
Using the elementary bound $\log (1+u)<u$ for $0<u<1,$ we see that
the LHS of (\ref{aritheq}) is $<\frac{\log p_{n+1}}{ k_n}.$
Combining the bounds on the LHS and the RHS we obtain
$$k_n \log k_n <p_{n+1} \log p_{n+1}.$$
Since the function $x \mapsto x \log x$ is non decreasing for $x
> >e$ we obtain $k_n <p_{n+1},$ that is
$$p_{n+1}-p_n>C\sqrt{p_n} \log_3 p_n,$$
which contradicts Cram\'er conjecture \cite{C}
$$p_{n+1}-p_n=O(\log^2 p_n).$$
 \qed

 But is also not increasing, as the next Proposition shows
 unconditionally.
{\prop \label{viol} The inequality $u_n<u_{n+1}$ is violated for
infinitely many $n$'s.}

\pf Suppose that  $u_n<u_{n+1}$  for $n$ big enough. Then for $n$
large enough we have

$$ u_n \le e^\gamma.$$

If RH is true that is a contradiction by \cite[Th. 2 (a)]{N}. If
RH is false that contradicts \cite[Th. 2 (b)]{N}. \qed

Thus $u_n$ is not a monotone sequence for $n$ big enough.
\section{Analogous problem for Dedekind $\Psi$ function}
Recall that the Dedekind $\Psi$ function is the multiplicative
function defined by $$\Psi(n)=n\prod_{p \vert n}(1+\frac{1}{p}).$$
Define the sequence $v_n=\frac{\Psi(N_n)}{N_n \log_2 N_n}.$ We
proved in \cite{SP} the two statements
\begin{itemize}
\item $v_n> \frac{e^\gamma}{\zeta(2)}$  for all $n\ge 3$ iff RH is
true \item $\lim v_n=\frac{e^\gamma}{\zeta(2)}$
\end{itemize}

Thus, like for the sequence $u_n$ it is natural to wonder if $v_n$
is decreasing.

{\prop The inequality $u_n>u_{n+1}$ is equivalent to
\begin{equation} \label{arith2}
\log(1+\frac{\log p_{n+1} }{\theta(p_n)}
                                        )>\frac{\log \theta(p_{n})}{p_{n+1}}
\end{equation}
}

\pf The inequality $v_n>v_{n+1}$ can be written as

$$\frac{\Psi(N_{n})}{N_n\log_2 N_n}  > \frac{\Psi(N_{n+1})}{N_{n+1}\log_2
N_{n+1}}.$$

 Note first that
$$ \frac{\Psi(N_{n+1})}{N_{n+1}}=(1+1/p_{n+1})\frac{\Psi(N_n)}{N_{n}},$$

so that, after clearing denominators, $v_n>v_{n+1}$ is equivalent
to
$$ \log_2( N_{n+1})> \log_2 N_n (1+1/p_{n+1}), $$
or, distributing, to
$$ \log_2(N_{n+1} )- \log_2 N_n> \frac{\log_2 N_{n}}{p_{n+1}}. $$
Like in the proof of Proposition \label{arith} we have
$$\log_2(N_{n+1} )=\log_2 N_n+ \log(1+\frac{\log p_{n+1}}{\log N_n}).$$
Combining the last two statements we obtain
$$ \log(1+\frac{\log p_{n+1}}{\log N_n} )>\frac{\log_2 N_{n}}{p_{n+1}}.$$
 The result follows then upon letting $\log N_n=\theta(p_n).$
\qed

Note that inequality \ref{arith2} is slightly looser than
inequality \ref{aritheq}. Still,
 the analogue of Proposition \ref{cram} is true:

{\prop \label{cram2} The inequality $v_n>v_{n+1}$ is violated for
infinitely many $n$'s.}

Similarly one can prove the analogue of Proposition \ref{viol} by
using the arguments in the proof of
\cite[Th. 2]{SP}.

 {\prop \label{viol2} The inequality
$v_n<v_{n+1}$ is violated for infinitely many $n$'s.}

The proofs of Propositions \ref{cram2} and \ref{viol2} are
completely analogous to the case of Euler $\varphi$ and are
omitted.

\section*{Acknowledgements}The third author acknowledges the hospitality of Postech Math Dept where this work was
performed.

\bibliographystyle{amsplain}

\end{document}